\documentclass[12pt]{smfart}

\usepackage{smfthm}
\usepackage{smfenum}
\usepackage{amssymb}
\usepackage{amscd}

\newcommand{\Qp}{\mathbf{Q}_p}
\newcommand{\Zp}{\mathbf{Z}_p}
\newcommand{\Cp}{\mathbf{C}_p}

\newcommand{\ra}{\rightarrow}

\newcommand{\Qpbar}{\overline{\mathbf{Q}}_p}
\newcommand{\Qpnrhat}{\widehat{\mathbf{Q}}_p^{\mathrm{nr}}}

\newcommand{\on}{\operatorname}
\newcommand{\OO}{\mathcal{O}}

\renewcommand{\hat}{\widehat}
\renewcommand{\tilde}{\widetilde}
\renewcommand{\phi}{\varphi}

\newcommand{\ZZ}{\mathbf{Z}}

\newcommand{\bcris}{\mathbf{B}_{\mathrm{cris}}}

\newcommand{\bdR}{\mathbf{B}_{\mathrm{dR}}}

\newcommand{\bplus}{\mathbf{B}^+}

\newcommand{\aplus}{\mathbf{A}^+}

\newcommand{\bhol}[1]{\mathbf{B}^+_{\mathrm{rig} #1}}

\newcommand{\dcris}{\mathbf{D}_{\mathrm{cris}}}
\newcommand{\ddR}{\mathbf{D}_{\mathrm{dR}}}

\newcommand{\repcris}[1]{\mathbf{Rep}_{\mathrm{cris}}^{#1}}
\newcommand{\modfil}[1]{\mathbf{MF}_{\OO_F}^{#1}}
\newcommand{\fil}{\mathrm{Fil}^0}

\title{Nombres de Tamagawa de certaines repr\'esentations cristallines}

\author{Laurent Berger}

\address{Harvard Dept of Mathematics  \\
      One Oxford Street \\
      Cambridge, MA 02138-2901 \\ USA}
\email{laurent@math.harvard.edu}
\urladdr{www.math.harvard.edu/\~{}laurent}


\thanks{This research was partially conducted by the author for the Clay
  Mathematics Institute.}

\date{Septembre 2002}
\begin{document}

\begin{abstract}
L'objet de ce texte est de montrer la conjecture $C_{EP,F}(V)$ pour
certaines repr\'esentations semi-stables, et la conjecture $\delta_{\Zp}(V)$
pour certaines repr\'esentations cristallines. Les deux ingr\'edients
principaux sont d'une part le calcul de nombres de Tamagawa en
utilisant une variante des r\'esultats de
Fontaine-Laffaille et de Bloch-Kato, et d'autre part les r\'esultats
de continuit\'e de l'exponentielle de Bloch-Kato montr\'es par Perrin-Riou.
\end{abstract}

\begin{altabstract}
The purpose of this article is to give a proof of the $C_{EP,F}(V)$
conjecture for some semi-stable representations and of the
$\delta_{\Zp}(V)$ conjecture for some crystalline
representations. There are two major ingredients: first, the
computation of Tamagawa numbers using a variation on
results of Fontaine-Laffaille
and of Bloch-Kato, and second some continuity results for Bloch-Kato's
exponential which were proved by Perrin-Riou.
\end{altabstract}

\maketitle
\setcounter{tocdepth}{2}
\tableofcontents

\setlength{\baselineskip}{18pt}

\section*{Introduction}
Soit $p$ un nombre premier $\neq 2$ et 
$F$ une extension finie non ramifi\'ee de $\Qp$.
L'objet de ce texte est de montrer
les conjectures $C_{EP,F}(V)$ de Fontaine et Perrin-Riou
pour les repr\'esentations
semi-stables $V$, dont tous les sous-quotients irr\'eductibles $W$
satisfont les propri\'et\'es suivantes:
\begin{enumerate}
\item $W$ est cristalline;
\item la longueur de la filtration sur $\ddR(W)$ est $\leq p-1$;
\item $W$ satisfait au moins une des deux conditions ci-dessous:
\begin{enumerate}
\item les poids de Hodge-Tate de $W$ sont dans $[-(p-2);p-1]$; 
\item $\phi : \dcris(W) \ra \dcris(W)$ est semi-simple
\footnote{Si $E$ est un $K$-espace vectoriel,
on dit que $f: E \ra E$ est semi-simple 
en $\alpha \in K$ si $\on{ker} (f-\alpha) \cap
  \on{im} (f-\alpha) = 0$.} en $p^{-j}$
  pour tout $j \in \ZZ$ (on dira dans la suite 
que $W$ est $\phi$-semi-simple).
\end{enumerate}
\end{enumerate}
Comme corollaire et ingr\'edient de la d\'emontration, on montre la
conjecture $\delta_{\Zp}(V)$ de Perrin-Riou pour les repr\'esentations
cristallines $V$, dont tous les sous-quotients irr\'eductibles $W$
satisfont les propri\'et\'es suivantes:
\begin{enumerate}
\item $W$ est cristalline;
\item la longueur de la filtration sur $\ddR(W)$ est $\leq p-1$; 
\item $W$ est $\phi$-semi-simple.
\end{enumerate}

Rappelons le contenu de la conjecture $C_{EP,F}(V)$. 
Si $V$ est une repr\'esentation 
semi-stable de $G_F$, soient $\dcris(V)$, $\ddR(V)$ et $t(V) = \ddR(V)
/ \on{Fil}^0 \ddR(V)$ les invariants associ\'es \`a $V$ par la
th\'eorie de Fontaine, via les anneaux $\bcris$ et $\bdR$. 
Dans la suite, $T$ d\'enotera un $\Zp$-r\'eseau
de $V$. Comme $t(V^*(1))^* \simeq \on{Fil}^0 \ddR(V)$,
on a une suite exacte:
\[ \begin{CD}
0 @>>> t(V^*(1))^* @>>> \ddR(V) @>>> t(V) @>>> 0 
\end{CD} \]
qui permet d'identifier ${\det}_{\Qp} \ddR(V)$ \`a ${\det}_{\Qp}t(V) \otimes
{\det}_{\Qp}^{-1} t(V^*(1))$, et une suite exacte (cf. \ref{ebks} (\ref{sfv}))
\begin{multline*}
\begin{CD}
0 @>>> H^0(F,V) @>>> \dcris(V) @>{1-\phi,\lambda_V}>> \dcris(V) \oplus
t(V) @>{\exp_{F,V}}>> H^1(F,V)   @>{\exp^*_{F,V^*(1)}}>>
\end{CD} \\
\begin{CD}  \dcris(V^*(1))^* \oplus t(V^*(1))^*
  @>{1-{}^t\phi^{-1}, {}^t \lambda_{V^*(1)}}>> \dcris(V^*(1))^* @>>>
  H^2(F,V) @>>> 0
\end{CD}
\end{multline*}
provenant de la suite exacte fondamentale et de sa duale, qui
d\'efinit des isomorphismes $\tilde{e}_V : \Delta_F(V) \simeq
\tilde{\Delta}_{EP,F}(V)$, o\`u
\[  \Delta_F(V) = {\det}^{-1}_{\Qp} \ddR(V) \otimes 
{\det}_{\Qp} \on{Ind}_{\Qp}^F V,\]
et
\[ \tilde{\Delta}_{EP,F}(V) = \otimes_{0 \leq i \leq 2}
   {\det}_{\Qp}^{(-1)^i} H^i(F,V) \otimes
{\det}_{\Qp} \on{Ind}_{\Qp}^F V. \]
Soit
\[ \Delta_{EP,F,\Zp}(V) = \tilde{e}_V^{-1} \left( \otimes_{0 \leq i \leq 2}
   {\det}_{\Qp}^{(-1)^i} H^i(F,T) \otimes
{\det}_{\Qp} \on{Ind}_{\Qp}^F T \right), \]
c'est un sous $\Zp$-module de rang $1$ de $\Delta_F(V)$, et la
conjecture $C_{EP,F}(V)$ propose une formule qui permet de le calculer. Posons
$h_j(V) = \dim_F (\on{Fil}^j \ddR(V) / \on{Fil}^{j+1} \ddR(V))$, 
$t_H(V) = \sum_{j \in \ZZ} j \cdot h_j(V)$, et
\[ \Gamma^*(j) =
\begin{cases}
(j-1)! & \text{si $j \geq 1$,} \\
(-1)^j / (-j)! & \text{si $j \leq 0$.}
\end{cases} \]
Si $\omega \in \Delta_F(V)$, soient $\tilde{\xi}_V(\omega) \in \Qpbar$ 
le coefficient de $t^{-t_H(V)}$
dans l'image
de $\omega$ par l'isomor\-phisme de comparaison $\Delta_F(V) \ra
t^{-t_H(V)} \Qpbar$ et $\eta_V(\omega)=
|\tilde{\xi}_V(\omega)|_p$
(si $V$ est semi-stable, la d\'efinition de $\eta_V$ est plus
compliqu\'ee, mais nous n'utiliserons que le cas cristallin).

La conjecture $C_{EP,F}(V)$ de Fontaine et Perrin-Riou
peut alors s'\'enoncer ainsi:
\begin{conj}[$C_{EP,F}(V)$]
Pour tout $\omega \in \Delta_F(V)$, on a
\[ \Delta_{EP,F,\Zp}(V) = \Zp \cdot \det \left(-\phi \mid
\dcris(V^*(1)) \right) \prod_{j \in \ZZ} \Gamma^*(-j)^{-h_j(V)[F:\Qp]}
\eta_V(\omega) \omega. \]
\end{conj}
Le r\'esultat principal de cet article est que la conjecture ci-dessus
est vraie pour une repr\'esentation cristalline irr\'eductible, dont
la longueur de la filtration est $\leq p-1$, et qui est
soit $\phi$-semi-simple, soit \`a poids de Hodge-Tate dans $[-(p-2);p-1]$. 
Comme la conjecture $C_{EP,F}(V)$
est stable par suite exactes, on en d\'eduit:
\begin{theo}
La conjecture $C_{EP,F}(V)$ est vraie
pour les repr\'esentations 
semi-stables $V$ de $G_F$, telles que tous 
les sous-quotients irr\'eductibles de
$V$ sont des repr\'esentations cristallines dont la 
longueur de la filtration est $\leq p-1$, et qui sont
soit $\phi$-semi-simples, soit \`a poids de Hodge-Tate dans $[-(p-2);p-1]$.
\end{theo} 

Comme ingr\'edient et corollaire de nos calculs, nous montrons aussi
la conjecture $\delta_{\Zp}(V)$ de Perrin-Riou 
(on renvoie \`a \ref{deltazp} pour des rappels sur les constructions
de Perrin-Riou) pour les 
repr\'esentations cristallines irr\'eductibles $\phi$-semi-simples dont
la longueur de la filtration est $\leq p-1$. 
Comme la conjecture $\delta_{\Zp}(V)$
est stable par suite exactes, on en d\'eduit:

\begin{theo}
La conjecture $\delta_{\Zp}(V)$ est vraie
pour les repr\'esentations 
cristallines $V$ de $G_F$, telles que tous les sous-quotients 
irr\'eductibles de
$V$ sont des repr\'esenta\-tions cristallines $\phi$-semi-simples dont la 
longueur de la filtration est $\leq p-1$.
\end{theo}

Remarquons que l'on conjecture que les repr\'esentations qui
``viennent de la g\'eom\'etrie'' sont $\phi$-semi-simples.
Disons quelques mots sur la d\'emonstration. Ces r\'esultats sont
d\'ej\`a connus pour les repr\'esentations $V$ qui n'ont qu'un poids
de Hodge-Tate,
et on suppose donc que $V$ n'a pas de sous-quotient
qui est une tordue d'une repr\'esentation non-ramifi\'ee. 
On commence par calculer les nombres de Tamagawa de $V$, pour de telles
repr\'esentations $V$ satisfaisant de nombreuses restrictions sur
leurs poids de Hodge-Tate. On en d\'eduit certains cas de la
conjecture $C_{EP,F}$. Ensuite, on utilise la construction de
Perrin-Riou pour d\'eterminer comment se comporte $C_{EP,F}(V)$ quand
on twiste $V$. Ceci nous permet de d\'emontrer simultan\'ement
$C_{EP,F}(V)$ et $\delta_{\Zp}(V)$.


\renewcommand{\thesection}{\Roman{section}}

\section{Repr\'esentations cristallines}

\Subsection{Rappels et notations}
On se donne un nombre premier $p \neq 2$, et $F$ est une extension
finie non-ramifi\'ee de $\Qp$. Soit
$\Qpbar$ une
cl{\^o}ture alg{\'e}brique de $\Qp$ et 
$\Cp=\hat{\Qpbar}$ sa compl{\'e}tion
$p$-adique. 
On pose $G_F=\on{Gal}(\Qpbar/F)$.
On pose aussi $F_n=F(\mu _{p^n})$ et 
$F_{\infty}$ est d{\'e}fini comme 
{\'e}tant la r{\'e}union des $F_n$.  
Soit $H_F$ le noyau du caract{\`e}re cyclotomique 
$\chi: G_F \ra \Zp^*$ et $\Gamma_F=G_F/H_F$
le groupe de Galois de $F_{\infty}/F$
qui s'identifie via le caract{\`e}re cyclotomique {\`a} $\Zp^*$.

Dans la suite, $V$ sera une repr\'esentation $p$-adique de
$G_F$. Soient $\bcris$ et $\bdR$ les anneaux de p\'eriodes
construits par Fontaine. Rappelons que $\bcris$ est une $F$-alg\`ebre
munie d'un Frobenius $\phi$ et que $\bdR$ est une $F$-alg\`ebre
munie d'une filtration. Les principales propri\'et\'es dont nous
ferons usage sont le fait que l'on a une inclusion
$\lambda: \bcris \ra \bdR$, et la ``suite
exacte fondamentale'': 
\[ \begin{CD} 0 @>>> \Qp @>>> \bcris @>{1-\phi,\lambda}>> \bcris
  \oplus \bdR/\bdR^+ @>>> 0. \end{CD} \] 
Si $V$ est une repr\'esentation $p$-adique,
on pose $\dcris(V) = (\bcris \otimes_{\Qp} V)^{G_F}$, c'est un
$F$-espace vectoriel de dimension $\leq d = \dim_{\Qp}(V)$ et 
on dit que $V$ est cristalline s'il y a \'egalit\'e. On d\'efinit de
m\^eme $\ddR(V)$ et les repr\'esentations de de Rham. Le $F$-espace
vectoriel $\dcris(V)$ est muni d'un Frobenius $\phi$, et $\ddR(V)$ est
muni d'une filtration $\on{Fil}$ (qui d\'efinit une filtration sur
$\dcris(V)$ via l'inclusion $\lambda_V : \dcris(V) \ra \ddR(V)$). Les
entiers $j \in \ZZ$ tels que $\on{Fil}^{-j} 
\ddR(V) \neq \on{Fil}^{-j+1} \ddR(V)$
sont appel\'es les poids de Hodge-Tate de $V$. 

On dira que $V$ est $\phi$-semi-simple si $\phi : \dcris(V) \ra
\dcris(V)$ est semi-simple en $p^{-j}$ pour tout $j \in \ZZ$.

Dans tout cet article, $a$ et $b$ seront deux nombres entiers tels que
$a \leq 0$,  $1 \leq b$ et $b-a \leq p-1$. On s'int\'eressera aux
repr\'esentations cristallines dont les poids de Hodge-Tate sont dans 
$[a;b]$.

\Subsection{Construction de r\'eseaux}
Nous aurons besoin de rappeler et pr\'eciser certaines des
constructions de \cite{LBcr}. 
Les r\'esultats que nous obtenons sont de l\'eg\`eres variantes de
``Fontaine-Laffaille'' (cf. \cite{FL82}).

Pour $a \leq b \in \ZZ$,
soit $\modfil{[a;b]}$ la cat\'egorie (ab\'elienne) des 
$\OO_F$-modules $D$ libres de type fini, munis d'une
filtration d\'ecroissante $\{\on{Fil}^i D\}_i$
telle que $\on{Fil}^{-b}
D =D$ et $\on{Fil}^{-a+1} D = 0$, et d'une application 
$\phi : D \ra D$ telle que $D = \sum_{i \in \ZZ} p^{-i} \phi(\on{Fil}^i
D)$. On dira que $D$ est fortement divisible.
Soit $\modfil{[a;b*]}$ la sous-cat\'egorie pleine de 
$\modfil{[a;b]}$ constitu\'ee des objets $D$ tels que $\phi :D \ra D$ 
est de pente $>-b$ (il est \'equivalent de demander que 
le $\phi$-module sous-jacent \`a $D$ n'a pas de
partie de pente $-b$). De m\^eme, 
$\modfil{[a*;b]}$ est la sous-cat\'egorie pleine de 
$\modfil{[a;b]}$ constitu\'ee des objets $D$ tels que $\phi :D \ra D$ 
est de pente $<-a$ et $\modfil{[a*;b*]} = 
\modfil{[a*;b]} \cap \modfil{[a;b*]}$.

Ensuite, on d\'efinit $\repcris{[a;b]}$ comme \'etant 
la cat\'egorie des $\Zp$-modules libres de type fini $T$ qui sont
des r\'eseaux d'une repr\'esentation cristalline $V$ \`a poids de
Hodge-Tate dans $[a;b]$. Soit $\repcris{[a;b*]}$ la sous-cat\'egorie
pleine de $\repcris{[a;b]}$ constitu\'ee des objets $T$ tels que $V$
n'a pas de sous-objet isomorphe \`a $V_0(b)$ avec $V_0$
non-ramifi\'ee. De m\^eme, $\repcris{[a*;b]}$ est la sous-cat\'egorie
pleine de $\repcris{[a;b]}$ constitu\'ee des objets $T$ tels que $V$
n'a pas de quotient isomorphe \`a $V_0(a)$ avec $V_0$ non-ramifi\'ee,
et $\repcris{[a*;b*]} = \repcris{[a*;b]} \cap \repcris{[a;b*]}$.

Si $\aplus_F=\OO_F[[\pi]]$, muni d'un Frobenius $\phi$ et d'une action
de $\Gamma_F$ semi-lin\'eaires et d\'efinis par les formules
\[ \phi(\pi)=(1+\pi)^p-1 \quad\text{et}\quad
\gamma(\pi)=(1+\pi)^{\chi(\gamma)}-1, \]
rappelons que l'on a construit dans \cite{LBcr}
un foncteur $T \mapsto N(T)$, de 
$\repcris{[a;b]}$ dans la cat\'egorie des $\aplus_F$-modules libres de
type fini, munis d'une action de $\Gamma_F$ triviale sur $N(T)/\pi$ et
d'une application $\phi : \pi^b N(T) \ra \pi^b N(T)$ telle que
$\pi^b N(T) / \phi^*(\pi^b N(T))$ est tu\'e par $q^{b-a}$
(ce foncteur est d\'efini dans \cite{LBcr} pour les repr\'esentations 
$V = \Qp \otimes_{\Zp} T$ dont les poids de Hodge-Tate sont
n\'egatifs, et on l'\'etend par la formule $N(T)=\pi^{-h} N(T(h))$).

On a d\'efini pour $i\in\ZZ$, 
$\on{Fil}^i N(T) = \{ x \in N(T), \phi(x) \in q^i N(T)
\}$. Revenons au cas o\`u $b-a \leq p-1$.
\begin{lemm}
Le $\phi$-module $D = N(T)/\pi$, muni de la filtration induite, est un
objet de $\modfil{[a;b]}$.
\end{lemm}

\begin{proof}
C'est un r\'esultat de Wach (cf \cite[th\'eor\`eme
3]{Wa97}). Donnons-en une d\'emon\-stration 
r\'edig\'ee un peu diff\'eremment. 
Quitte \`a tordre, on se ram\`ene au cas o\`u $b=0$. Il
est trivial que $\phi(\on{Fil}^i D) \subset p^i D$, et il suffit donc
de voir que la filtration est admissible ou, ce qui revient au m\^eme,
que la filtration $\on{Fil}_V$
induite par celle de $N(V)/\pi$ co\"{\i}ncide avec
celle de $N(T)/\pi$. Nous allons donc montrer que $\on{Fil}_V^i D
\subset \on{Fil}^i D$ (l'autre inclusion est triviale), pour
$1 \leq i \leq p-1$ (si $i \geq p$, alors $\on{Fil}_V^i D = 0$).

Si $\gamma$ est un g\'en\'erateur de $\Gamma_F$,
posons \[ T_i(\gamma)=(1-\chi(\gamma)^{-1}\gamma)
(1-\chi(\gamma)^{-2}\gamma) \cdots (1-\chi(\gamma)^{-(i-1)}\gamma). \]
Si $\overline{x} \in N(T)$ est dans $\on{Fil}_V^i N(T)$, c'est donc
qu'il  est l'image modulo $\pi$ de 
$x \in N(V)$ tel que $\phi(x) \in q^i N(V)$, et tel que
l'on puisse \'ecrire $x =x_0 + \pi y_1$ avec $x_0 \in N(T)$ et $y_1
\in N(V)$. On a $T_i(\gamma)x = T_i(\gamma) x_0 + \pi^i y_i$ 
avec $y_i \in N(V)$. Comme $T_i(\gamma)x \in \on{Fil}^i N(V)$, 
et comme $\pi^i y_i \in \on{Fil}^i N(V)$ trivialement, on a
$T_i(\gamma)x_0 \in \on{Fil}^i N(V) \cap N(T) = \on{Fil}^i N(T)$, ce
qui fait que $T_i(\gamma) \overline{x} \in \on{Fil}^i D$. 
Si $i \leq p-1$, alors $T_i(\gamma)$ agit par une unit\'e $p$-adique
sur $D$ et donc on a bien $x \in \on{Fil}^i D$. 
\end{proof}

\begin{lemm}\label{isomcomp}
Si $T\in \repcris{[a;b]}$ 
est un r\'eseau d'une repr\'esentation cristal\-line 
$V$ et $D=N(T)/\pi$,
alors le d\'eterminant de l'isomorphisme de comparaison
$\bcris \otimes_{\Zp} T \simeq \bcris \otimes_{\OO_F} D$ appartient
\`a $t^{t_H(V)} \OO_{\Qpnrhat}^*$.
\end{lemm}

\begin{proof}
Tout d'abord, par \cite{LBcr},
le d\'eterminant de l'inclusion $\aplus 
\otimes_{\aplus_F} N(T) \subset \aplus \otimes_{\Zp} T$ appartient \`a
$\pi^{t_H(V)} (\aplus)^*$, et donc \`a $\pi^{t_H(V)} (\aplus_F)^*
 \OO_{\Qpnrhat}^*$.

Ensuite, rappelons que 
$\dcris(V)= (\bhol{,F} \otimes_{\aplus_F}
N(T))^{\Gamma_F}$, et donc que $D \subset \dcris(V)$ s'identifie
\`a l'ensemble des $x \in (\bhol{,F} \otimes_{\aplus_F}
N(T))^{\Gamma_F}$ tels que $x(0) \in N(T) / \pi$. 
Le d\'eterminant de l'inclusion
$\bhol{,F} \otimes_{\OO_F} D \subset \bhol{,F} \otimes_{\aplus_F}
N(T)$ calcul\'e dans des bases de $D$ et $N(T)$ 
est donc \'egal \`a $(\log(1+\pi)/\pi)^{t_H(V)} d(\pi)$ 
o\`u $d(\pi) \in (\bplus_F)^*$ et 
$d(0) \in \OO_F^*$. On en d\'eduit le r\'esultat.
\end{proof}

\begin{prop}\label{constfl}
Si $D \in \modfil{[a;b*]}$, alors il existe $T \in \repcris{[a;b]}$ tel
que $N(T)/\pi \simeq D$. 
\end{prop}

\begin{proof}
Commen\c{c}ons par remarquer que si 
$\mu=(p/(q-\pi^{p-1})) \in \aplus_F$, alors 
$\mu$ est inversible dans $\aplus_F$, et pour tout $s \geq 0$,
$\mu^s q^s = p^s \mod{\pi^{p-1}}$.

Quitte \`a tordre $D$, on se ram\`ene au cas o\`u 
$b=0$ et donc
$\on{Fil}^0 D = D$. 
Si $0 \leq r_1 \leq \cdots \leq r_d \leq p-1$ sont les sauts de la
filtration, alors dans une base de $D$ adapt\'ee \`a
la filtration, la matrice de $\phi$ est  
\'egale \`a $P_0 A$ o\`u $P_0=\on{Diag}(p^{r_1},\cdots,p^{r_d})$ et 
$A \in \on{GL}(d,\OO_F)$. On pose
$P=\on{Diag}(q^{r_1}\mu^{r_1},\cdots,q^{r_d}\mu^{r_d})A$ 
ce qui fait que $P$ est une matrice \`a coefficients dans $\aplus_F$, 
\'egale modulo $\pi^{p-1}$ \`a la matrice de $\phi$ sur $D$. 
De plus, si $\gamma \in \Gamma_F$, alors $\gamma(P^{-1})P
\in \on{Id}+\pi^{p-1} \on{M}(d,\aplus_F)$.  On \'ecrira $\gamma(P^{-1})P
= \on{Id} + \pi^{p-1} Q$.

Pour terminer la preuve, il suffira de construire 
pour $\gamma \in \Gamma_F$ une matrice $H = H_\gamma \in
\on{M}(d,\aplus_F)$ telle que l'on puisse d\'efinir l'action de
$\gamma$ sur $N(T)$ par $G_\gamma=\on{Id}+\pi^{p-1}H_\gamma$, c'est-\`a-dire
telle que $\gamma(P)G_\gamma=\phi(G_\gamma)P$. Cette \'equation est
\'equivalente \`a:
\[ H-q^{p-1}\gamma(P^{-1}) \phi(H) P =
\frac{\gamma(P^{-1})P-\on{Id}}{\pi^{p-1}} = Q. \]
Il suffira donc de montrer que l'application 
$X \mapsto X-q^{p-1}\gamma(P^{-1}) \phi(X) P$ 
de $\on{M}(d,\aplus_F)$ dans lui-m\^eme est surjective. Il est clair
que si $Y \in \pi \on{M}(d,\aplus_F)$, alors la s\'erie
\[ Y+ q^{p-1}\gamma(P^{-1}) \phi(Y) P
+ q^{p-1}\gamma(P^{-1}) \phi(q^{p-1}\gamma(P^{-1})) \phi^2(Y) \phi(P)
P + \cdots \]
converge vers $X \in \pi \on{M}(d,\aplus_F)$ tel que
$X-q^{p-1}\gamma(P^{-1}) \phi(X) P = Y$. Il suffit donc de montrer que
l'application $X \mapsto X-p^{p-1}(P_0 A)^{-1} \phi(X) P_0 A$ 
de $\on{M}(d,\OO_F)$ dans lui-m\^eme est surjective. Le fait que 
$\phi: D \ra D$ n'a pas de partie de pente $0$ implique que 
$\prod_{i=0}^n \phi^{n-i}(P_0 A) \ra 0$ quand $n \ra \infty$ et
donc que si $Y \in \on{M}(d,\OO_F)$, alors la s\'erie
 \[ Y+ p^{p-1}(P_0 A)^{-1} \phi(Y) P_0 A
+ p^{p-1}(P_0 A)^{-1} \phi(p^{p-1}(P_0 A)^{-1}) \phi^2(Y) 
\phi(P_0 A) P_0 A + \cdots \]
converge vers $X \in \on{M}(d,\OO_F)$ tel que 
$X-p^{p-1}(P_0 A)^{-1} \phi(X) P_0 A = Y$.
\end{proof}

\begin{rema}
La m\^eme d\'emonstration montre que 
si $D \in \modfil{[a*;b]}$, alors il existe $T \in \repcris{[a;b]}$ tel
que $N(T)/\pi \simeq D$. En effet, dans la d\'emonstration ci-dessus,
si $D$ n'a pas de partie de pente $p-1$, alors 
$\prod_{i=0}^n \phi^i(p^{p-1} (P_0 A)^{-1}) \ra 0$ quand $n \ra \infty$.
\end{rema}

\begin{coro}
Le foncteur $T \mapsto N(T)/\pi$ 
est une \'equivalence de cat\'egories
de $\repcris{[a;b*]}$ dans $\modfil{[a;b*]}$
et de $\repcris{[a*;b]}$ dans $\modfil{[a*;b]}$.
\end{coro}

\begin{proof}
La proposition \ref{constfl} ci-dessus en construit un quasi-inverse.
\end{proof}

\begin{defi}
On \'ecrira abusivement $\dcris(T) = N(T) / \pi$.
\end{defi}

\begin{coro}\label{fl}
Si $T \in \repcris{[a;b*]}$ 
est un r\'eseau d'une repr\'esentation cristal\-line 
$V$, alors
on a des isomorphismes canoniques \[ \dcris(T)/(1-\phi) \fil \dcris(T) 
\simeq \on{Ext}^1_{\modfil{[a;b]}}(\on{Id}_D,\dcris(T)) 
\simeq H^1_f(F,T). \]
\end{coro}

\begin{proof}
Rappelons que $0 \in [a;b]$ ce qui fait que l'objet trivial $\on{Id}_D$ 
est dans $\modfil{[a;b]}$ et que $b \geq 1$, ce qui fait que
$\on{Id}_D \in \modfil{[a;b*]}$.
Le fait que  $T \in \repcris{[a;b*]}$ 
implique que $D(T)$ est aussi dans $\modfil{[a;b*]}$.
Le corollaire r\'esulte alors de l'\'equivalence de cat\'egories ci-dessus.
\end{proof}

\section{Calculs de nombres de Tamagawa}

\Subsection{L'exponentielle de Bloch-Kato}\label{ebks}
Rappelons que les anneaux de Fontaine $\bcris$ et $\bdR$ sont reli\'es
par la ``suite exacte fondamentale'':
\[ \begin{CD} 0 @>>> \Qp @>>> \bcris @>{1-\phi,\lambda}>> \bcris
  \oplus \bdR/\bdR^+ @>>> 0, \end{CD} \] 
et que si l'on tensorise par $V$ et que l'on
prend la suite exacte longue de cohomologie associ\'ee \`a
$(\cdot)^{G_F}$, on trouve le d\'ebut de la suite exacte de
l'introduction:
\begin{equation}\label{sfv} 
\begin{CD}
0 @>>> H^0(F,V) @>>> \dcris(V) @>{1-\phi,\lambda_V}>> \dcris(V) \oplus
t(V) @>{\exp_{F,V}}>> H^1(F,V)  
\end{CD} \end{equation}
L'image de $\exp_{F,V}$ est 
\[ H^1_f(F,V) =  \on{ker} \left(H^1(F,V) \ra H^1(F,\bcris
\otimes_{\Qp} V) \right), \]
et la deuxi\`eme moiti\'e de la suite exacte de l'introduction
s'obtient en dualisant, et en utilisant un th\'eor\`eme de Bloch et
Kato qui dit que $H^1_f(F,V)$ et $H^1_f(F,V^*(1))$ sont orthogonaux.

Si $V$ est une repr\'esentation 
cristalline, alors
$\dcris(V)/(1-\phi) \fil \dcris(V) \simeq H^1_f(F,V)$ et
par \cite[lemma 4.5]{BK91}:
\begin{prop}\label{expbk}
L'application que l'on d\'eduit de $1-\phi$:
\[ \begin{CD} \dcris(V)/ \fil \dcris(V) 
@>{1-\phi}>>
\dcris(V)/(1-\phi) \fil \dcris(V) \simeq H^1_f(F,V) \end{CD} \]
co\"{\i}ncide avec l'exponentielle de Bloch-Kato
$\exp_{F,V}: \dcris(V)/ \fil \dcris(V) \ra H^1_f(F,V)$.
\end{prop}
Ceci nous donne une suite exacte
\begin{multline}\label{blochkato}
\begin{CD} 
0 @>>> H^0(F,V) @>>> \dcris(V)^{\phi=1} @>>> \dcris(V)/ \fil
\dcris(V) \end{CD} \\ \begin{CD}  
@>{1-\phi}>> \dcris(V)/(1-\phi) \fil \dcris(V) @>>>   
\dcris(V)/(1-\phi) \dcris(V) @>>> 0.
\end{CD} \end{multline}
On a aussi une suite exacte b\^ete:
\begin{multline}\label{bete} 
\begin{CD}
0 @>>> \dcris(V)^{\phi=1} @>>> \dcris(V) \end{CD} \\ \begin{CD}
@>{1-\phi}>>  \dcris(V) @>>> 
\dcris(V)/(1-\phi) \dcris(V) @>>> 0. 
\end{CD} \end{multline}

\Subsection{Nombres de Tamagawa}
On suppose toujours que $0,1 \in [a;b]$ et que $b-a \leq p-1$.
Rappelons que par d\'efinition, 
si $\omega_{t(V)}$ est une base de 
$\det_{\Zp} (\dcris(T)/ \fil \dcris(T))$, alors
$\on{Tam}^0_p(V)$ est le coefficient de $\omega_{t(V)}^{-1}$
dans l'image de
$\det_{\Zp} H^0(F,T) \otimes \det^{-1}_{\Zp} H^1_f(F,T)$ 
dans $\det^{-1}_{\Zp} (\dcris(T)/ \fil \dcris(T))$ par
la suite exacte
\begin{multline*}\label{tamag}
\begin{CD}
0 @>>> H^0(F,V) @>>> \dcris(V) @>{1-\phi,\lambda_V}>>
\end{CD} \\ \begin{CD} 
\dcris(V) \oplus \dcris(V)/ \fil
\dcris(V) @>>> H^1_f(F,V) @>>> 0.
\end{CD}
\end{multline*}
Si l'on combine le corollaire \ref{fl}, la proposition \ref{expbk}, 
et les suites exactes \ref{blochkato} et \ref{bete}, on trouve que:
\begin{prop}\label{tamagunit}
Si $T \in \repcris{[a;b*]}$ 
est un r\'eseau d'une repr\'esentation 
cristalline $V$, alors
$\on{Tam}^0_p(V) \in \ZZ_p^*$.
\end{prop}
Rappelons la proposition \cite[C.2.6]{BP95} qui fait le lien entre les
nombres de Tamagawa et la conjecture $C_{EP,F}$:
\begin{prop}
Si $\omega = \omega_{t(V^*(1))} \otimes 
\omega_{t(V)}^{-1} \in \Delta_F(V)$, alors
\[ \Delta_{EP,F,\Zp}(V) = \Zp \cdot \det \left(-\phi \mid
\dcris(V^*(1)) \right) \frac{\on{Tam}^0_p(V)}{\on{Tam}^0_p(V^*(1))}
\omega. \]
\end{prop}
La proposition ci-dessus et le lemme \ref{isomcomp} 
nous permettent de montrer $C_{EP,F}(V)$:
\begin{prop}
Si $T \in \repcris{[a*;b*]}$ est un r\'eseau d'une repr\'esentation 
cristalline $V$, alors la conjecture $C_{EP,F}(V)$ est vraie.
\end{prop}

\begin{proof}
Comme les poids de Hodge-Tate de $V$ sont tous dans $[-(p-2);p-1]$,
les facteurs $\Gamma^*(-j)$ qui interviennent
dans la d\'efinition de $\Delta_{EP,F,\Zp}(V)$ sont des unit\'es $p$-adiques.
La conjecture $C_{EP,F}(V)$ suit alors du fait que
$\on{Tam}^0_p(V)$ et $\on{Tam}^0_p(V^*(1)) \in \ZZ_p^*$, puisque si $V$ 
est cristalline \`a poids dans $[a;b] \subset [-(p-2);p-1]$, alors 
$V^*(1)$ est cristalline et  $T^*(1) \in \repcris{[(1-b)*;(1-a)*]}$
avec $[1-b;1-a] \subset [-(p-2);p-1]$.
\end{proof}

\begin{coro}
Si $T \in \repcris{[0*;p-1*]}$ est un r\'eseau d'une repr\'esentation 
cristalline $V$, 
alors la conjecture 
$C_{EP,F}(V(i))$ est vraie pour tout $i \in \{0,-1,\cdots,-(p-2)\}$.
\end{coro}

\section{L'exponentielle de Perrin-Riou}
\Subsection{Rappels sur l'exponentielle de Perrin-Riou}
Soit $\Delta_F$ le sous-groupe de torsion de $\Gamma_F$ et
$\Gamma_1=\chi^{-1}(1+p\Zp)$ ce qui fait que $\Gamma_F \simeq \Delta_F
\times \Gamma_1$. On pose $\Lambda=\Zp[[\Gamma_F]]$ et 
$\mathcal{H}(\Gamma_F) = \Qp[\Delta_F] \otimes_{\Qp}
\mathcal{H}(\Gamma_1)$ o\`u $\mathcal{H}(\Gamma_1)$ est l'ensemble des
$f(\gamma_1-1)$, avec $\gamma_1 \in \Gamma_1$, o\`u $f(T) \in \Qp[[T]]$
est un s\'erie formelle de rayon de convergence $\geq 1$.
 
Rappelons que le r\'esultat principal de Perrin-Riou dans \cite{BP94}
est la construction, pour une repr\'esentation cristalline $V$, 
d'une famille d'applications \[ \Omega_{V,h} : \mathcal{H}(\Gamma_F)
\otimes_{\Qp} \dcris(V) \ra \mathcal{H}(\Gamma_F)
\otimes_{\Lambda} H^1_{Iw}(F,V) / V^{H_F} \] dont la propri\'et\'e
principale est qu'elles interpolent les applications ``exponentielle
de Bloch-Kato''. Plus pr\'ecis\'ement, pour $h,j \gg 0$, le diagramme
\[ \begin{CD}
\left( \mathcal{H}(\Gamma_F) 
\otimes_{\Qp} \dcris(V(j)) \right)^{\tilde{\Delta}=0} 
@>{\Omega_{V(j),h}}>> \mathcal{H}(\Gamma_F)
\otimes_{\Lambda} H^1_{Iw}(F,V(j)) / V(j)^{H_F} \\
@V{\Xi_{n,V(j)}}VV @V{\on{pr}_{F_n,V(j)}}VV \\
F_n \otimes_F \dcris(V) @>{(h+j-1)! \exp_{F_n,V(j)}}>> H^1(F_n,V(j)) 
\end{CD} \]
est commutatif, o\`u $\tilde{\Delta}$ et 
$\Xi_{n,V}$ sont deux applications un peu
compliqu\'ees \`a d\'ecrire. Nous utiliserons le fait que 
$\on{Tw}_1 \Omega_{V,h} (\on{Tw}_1 \otimes
d_{-1})=-\Omega_{V(1),h+1}$ o\`u $d_{-1}$ est une base de $\dcris(\Zp(-1))$.

Perrin-Riou a propos\'e une conjecture relativement au d\'eterminant
de $\Omega_{V,h}$. Posons \begin{multline*} 
\delta(\Omega_V) = \prod_{j \geq 1-h} 
(\ell_{-j})^{-\dim_{\Qp} \on{Fil}^j \dcris(V)}  \\ \times
\Omega_{V,h} \left[ {\det}_{\Lambda} (\Lambda \otimes_{\Qp} \dcris(V)) 
\otimes_{\Lambda}
{\det}_{\Lambda}^{-1}H^1_{Iw}(F,V) \otimes_{\Lambda} 
{\det}_{\Lambda}H^2_{Iw}(F,V) \right].
\end{multline*}
Perrin-Riou a montr\'e que $\delta(\Omega_V)$ ne d\'epend pas
de $h$, et \cite[th\'eo 3.4.2]{BP94} que
$\delta(\Omega_V) \in \Qp \otimes_{\Zp}
\Lambda$; elle a de plus
conjectur\'e que $\delta(\Omega_V)$ est une unit\'e de 
$\Qp \otimes_{\Zp} \Lambda$. Ceci a \'et\'e montr\'e par Colmez 
\cite{Co98}. Colmez a en fait montr\'e la conjecture R\'ec($V$), qui
affirme que pour les accouplements naturels sur $H^1_{Iw}(F,V) 
\times H^1_{Iw}(F,V^*(1))$ et sur 
$\mathcal{H}(\Gamma) \otimes_{\Qp} \dcris(V) 
\times  \mathcal{H}(\Gamma)
\otimes_{\Qp} \dcris(V^*(1))$, l'adjoint de $\Omega_{V,h}$
est essentiellement l'inverse de $\Omega_{V^*(1),1-h}$. 
Cela implique que $\delta(\Omega_V) \delta(\Omega_{V^*(1)})^\iota$
est une unit\'e de 
$\Qp \otimes_{\Zp} \Lambda$, et donc de m\^eme pour $\delta(\Omega_V)$
et $\delta(\Omega_{V^*(1)})$.

\Subsection{La conjecture $\delta_{\Zp}(V)$}\label{deltazp}
Soient $T$ un r\'eseau de $V$ et 
$M$ un r\'eseau de $\dcris(V)$, 
tels que le d\'eterminant de l'isomor\-phisme de comparaison 
entre $\bcris \otimes_{\Zp} T$ et $\bcris \otimes_{\OO_F} M$
appartient \`a $t^{t_H(V)} \OO_{\Qpnrhat}^*$,
et $\omega$ la base de $\Delta_F(V)$ que l'on en d\'eduit. 
Posons: \begin{multline*} 
\delta_{\Zp}(\Omega_V) = \prod_{j \geq 1-h} 
(\ell_{-j})^{-\dim_{\Qp} \on{Fil}^j \dcris(V)}  \\ \times
\Omega_{V,h} \left[ {\det}_{\Lambda} (\Lambda \otimes_{\Zp} M) 
\otimes_{\Lambda}
{\det}_{\Lambda}^{-1}H^1_{Iw}(F,T) \otimes_{\Lambda} 
{\det}_{\Lambda} H^2_{Iw}(F,T) \right].
\end{multline*}
On peut alors \'enoncer la conjecture $\delta_{\Zp}(V)$:
\begin{conj}[$\delta_{\Zp}(V)$]
On a $\delta_{\Zp}(\Omega_V) \in \Lambda^*$.
\end{conj}
On voit que $\Delta_F^*$, le dual de $\Delta_F$, 
est compos\'e des $\varpi^i$, $0 \leq i \leq
p-2$ o\`u $\varpi$ est le caract\`ere de Teichm\"uller. Pour $i \in
\ZZ / (p-1)\ZZ$, on notera $e_i$ l'idempotent associ\'e \`a
$\varpi^i$. Le fait que $\on{Tw}_1 \Omega_{V,h} (\on{Tw}_1 \otimes
d_{-1})=-\Omega_{V(1),h+1}$ implique:
\begin{lemm}
On a $\on{Tw}_1 \delta_{\Zp}(\Omega_V) = \delta_{\Zp}(\Omega_{V(1)})$.
En particulier, on a $e_{i+1} 
\delta_{\Zp}(\Omega_V) = 
e_i \delta_{\Zp}(\Omega_V(1))$.
\end{lemm}

L'ingr\'edient principal dont nous avons maintenant besoin
est une formule qui relie $\delta_{\Zp}(\Omega_V)$ et $C_{EP,F}(V)$, 
formule qui
est montr\'ee sous une forme un peu diff\'erente dans \cite{BP94}
(il faut faire attention au fait que les notations de \cite{BP94} et
de \cite{BP95} ne sont pas vraiment compatibles, la d\'efinition de
$\xi_V$ n'est par exemple pas la m\^eme):
\begin{prop}
Si $V$ est une repr\'esentation cristalline, 
et si $\phi$ est semi-simple en $1$ et $p^{-1}$ sur $\dcris(V)$, 
alors 
\[ \Delta_{EP,F,\Zp}(V) = \det( -\phi \mid \dcris(V^*(1))) 
 \prod_{j \in \ZZ} \Gamma^*(-j)^{-h_j(V)[F:\Qp]}
|\delta_{\Zp}(\Omega_V)(0)|_p \cdot \omega. \]
En particulier, la conjecture $e_0 \cdot 
\delta_{\Zp}(V)$ est \'equivalente \`a la
conjecture $C_{EP,F}(V)$.
\end{prop}

\begin{proof} 
Comme la conjecture R\'ec($V$) est maintenant d\'emontr\'ee, le
th\'eor\`eme \cite[3.5.4]{BP94} est vrai en tout g\'en\'eralit\'e (il
n'y a plus de ``grain de sable''). En prenant la $e_0$-composante de
la formule \cite[3.5.4]{BP94} 
de Perrin-Riou, on trouve la formule ci-dessus (si $f
\in \Qp \otimes_{\Zp} \Lambda$, 
alors $f(0)$ est par d\'efinition
la projection de $f$ dans $e_0 \cdot \Qp \otimes_{\Zp}
\Lambda / (1-\gamma_1) \simeq \Qp$). 
\end{proof}

On en d\'eduit donc que si $V$ est $\phi$-semi-simple, alors
la conjecture $\delta_{\Zp}(V)$
est \'equivalente \`a la r\'eunion des conjectures 
$C_{EP,F}(V(i))$ pour $i \in \{0,-1,\cdots,-(p-2)\}$. En particulier:
\begin{coro}
Si $T \in \repcris{[0*;p-1*]}$ est un r\'eseau d'une repr\'esentation 
cristalline $\phi$-semi-simple 
$V$, alors la conjecture 
$\delta_{\Zp}(V)$ est vraie.
\end{coro}
Le fait que  $\on{Tw_1}
\delta(\Omega_V) = \delta(\Omega_{V(1)})$ montre  que
$\delta_{\Zp}(V)$ et $\delta_{\Zp}(V(1))$ sont \'equivalentes. 
De plus, Perrin-Riou a montr\'e $\delta_{\Zp}(V)$ pour les tordues de
repr\'esentations non-ramifi\'ees. Ceci implique donc:
\begin{prop}
Si $V$ est une repr\'esentation cristalline $\phi$-semi-simple et
irr\'eductible, dont la longueur de la
filtration est $\leq p-1$, alors $\delta_{\Zp}(V)$ et $C_{EP,F}(V)$
sont vraies.
\end{prop}

Finalement, Perrin-Riou a montr\'e que les conjectures 
$\delta_{\Zp}$ (cf \cite[3.4.9]{BP94})
et $C_{EP,F}$ (cf \cite[C.2.9]{BP95})
sont stables par suites exactes, ce qui
fait que l'on a:

\begin{theo}
La conjecture $C_{EP,F}(V)$ est vraie
pour les repr\'esentations 
semi-stables $V$ de $G_F$, telles que 
tous les sous-quotients irr\'eductibles de
$V$ sont des repr\'esentations cristallines dont la 
longueur de la filtration est $\leq p-1$, et qui sont
soit $\phi$-semi-simples, soit \`a poids de Hodge-Tate dans $[-(p-2);p-1]$.
\end{theo} 

\begin{theo}
La conjecture $\delta_{\Zp}(V)$ est vraie
pour les repr\'esentations 
cristallines $V$ de $G_F$, telles que tous 
les sous-quotients irr\'eductibles de
$V$ sont des repr\'esenta\-tions cristallines $\phi$-semi-simples dont la 
longueur de la filtration est $\leq p-1$.
\end{theo}

Remarquons que 
l'on conjecture que si $V$ ``vient de la
g\'eom\'etrie'', alors $\phi$ est semi-simple sur $\dcris(V)$.

\end{document}